\documentclass[11pt,leqno]{article}
\usepackage{amsfonts, amsmath, amsthm, amssymb}
\newtheorem{theorem}{Theorem}[section]

\newtheorem{definitiontemp}[theorem]{Definition}

\newcommand{\half}{\displaystyle{\frac{1}{2}}}

\newcommand{\C}{{\mathbb C}}

\newcommand{\mf}[1]{\displaystyle{\mathfrak{#1}}}

\begin{document}

\title{ Center of Infinitesimal Hecke Algebra of $\mf{sl}_2$}
\author{Akaki Tikaradze      }
\maketitle

 In  this notes we are concerned with a special case of so called infinitesimal Hecke algebras introduced in [EGG] For simple Lie algebra $\mf{sp}(2n)$ they construct certain algebra (called Infinitesimal Hecke algebra) which is endowed with a filtration so that corresponding associated graded algebra is isomorphic to the enveloping algebta of  $\mf{U}(\mf{sp}(2n) \ltimes V)$, where $V$ is an obvious representation (this is a version of PBW theorem). Following is the exact definition for the easiest case of $\mf{sl}_2$.

  Let $\mathbb{C}$ be the field of complex numbers, and let $H_z$ be a semidirect product of $\mf{U}\mf{g}$ and $TV$ ($TV$ stands for tensor algebra), where by $\mf{g}$ we will denote $\mf{sl}_2$ and $V$ is its two dimensional representation over $\mathbb{C}$ with basis $x,y,$, mod out by the relation$ [x,y]=z,$ where $ z$ is an element of the center of $\mf{U}\mf{g}$. Thus $ H_z$ is an algebra depending on the parameter $ z,$ recall that $ [e,y]=x$, and $ x$ is a maximal vector( killed by $ e$) of weight 1 with respect to  $ h$ (when talking about maximal vectors and such we always consider adjoint $\mf{sl}_2$ action). Recall that one has an algebra filtration on $ H_z$ ($x,y$ have degree 1, and $ \mf{U}\mf{g}$ has degree 0 ), corresponding associated graded algebra is just  $\mf{U}(\mf{g} \ltimes V)$ which also is equal to  $H_0$, this  statement is a PBW theorem for this algebra (proved in [Kh]), and from now on we will denote this algebra just by $H$ (suppressing the index $0$). Representation theory of $H_z$ was extensively studied in [Kh]. Here we will be concerned with the center and derivations of this algebra.

   I would like to thank my advisor V. Ginzburg, and especially A. Khare for carefully going  through all the computations and enormous help with use of Tex. Partially suported by NSF grant DMS-0401164.

\begin{theorem}
Center of $H_z$ is a polynomial algebra of one variable, generating central element is of the filtration degree 2.

\end{theorem}

  \begin{proof}
First, recall an anti-isomrphism of $H _z$ which will be denoted by  $j$, defined as: $j(x)=y, j(y)=x, j(h)=h, j(e)=-f, j(f)=-e$ ([Kh]). lets note that this map fixes the following element $ey^2+hxy-fx^2$. Indeed, applying this map to this element we get 
\begin{eqnarray*}
-x^2f+xyh+y^2e=hxy+ey^2-fx^2-[e,y^2]-(-[f,x^2]) 
\end{eqnarray*}
but $[e,y^2]=xy+yx=[f,x^2]$ so this element is indeed fixed, also since $j$ is fixing $h$ and
center of $\mf{U}\mf{g}$, $j$ is fixing the algebra generated by the above elements, $ h$, and $\mf{Z}(\mf{U}\mf{g})$ (center of $\mf{U}L$). next step will be to exibit an element from this algebra which will 
commute with $e,x,h$
therefore will lie in the center of $H_z$. Let us denote the above element  by $t=ey^2+hxy-fx^2$, we have $[e,t]=e(xy+yx)-2exy+hx^2-hx^2=-ez$, so we see that $[e,t - \frac{1}{2} hz]=0$. We also have
\begin{eqnarray*}
[x,t] & = & e(zy+yz)-x^2y+hxz+yx^2\\
& = & 2ezy-e[z,y]+hzx-h[z,x]-(2zx-[z,x])\\
& = & 2ezy-2zx+[z,x]-e[z,y]+hzx-h[z,x],
 \end{eqnarray*} 
\begin{eqnarray*}
 [x,\half  hz]=-\half xz+\half h[x,z]=-\half zx+\half [z,x]+\half h[x,z]. 
\end{eqnarray*}
Also,
\begin{eqnarray*}
  [x,t-\half hz]=2ezy-\half 3zx+\half[z,x]-e[z,y]-\half h[z,x]+hzx.  
\end{eqnarray*} 
We denote the element $[x,t-\half hz]$ by $q_z$, now we want to produce an element $\omega_z$ in the center of $\mf{Z}(\mf{U}L)$ so that $[x,\omega_z]=q_z$, then $t-\half hz-\omega_z$ will be the desired central element. Let $\Delta=h^2+4ef-2h$ be a multiple of the Casimir element, next we will  compute  commutators of powers of $\Delta$ with $x$ and $y$, we have
\[ [\Delta,x]=hx+xh+4ey-2x=(2h-3)x+4ey \]
\[ [\Delta,y]=-hy-yh+4xf+2y=-(2hy+y)+4(fx-y)+2y=(-2h-3)y+4fx \].

Let us put $ [\Delta^n,x]=(f_n(\Delta)h+g_n(\Delta))x+2f_n(\Delta)ey $, ($f_n, g_n$ are yet undetermined polynomials). It is always possible to write them in this way because $[\Delta^n,x]$ is a maximal vector of weight 1, so coefficient in fromt of $y$ is maximal vector of weight 2, by classification of maximal vectors in $\mf{U}\mf{g}$ (explained below), this element is a product of $e$ and an element of a center, similarly we treat coefficient in front of $x$. We have

\begin{align*}
[\Delta^{n+1},x] =& \Delta^n[\Delta,x]+[\Delta^n,x]\Delta\\
 = & \Delta^n(2h-3)x+\Delta^n4ey+(f_n(\Delta)h+g_n(\Delta))x\Delta +2f_n(\Delta)ey\Delta\\
 = & (\Delta^n(2h-3)+f_n(\Delta)\Delta h+g_n(\Delta)\Delta)x+(2f_n(\Delta)\Delta e+\Delta^n4e)y\\
 &  - (f_n(\Delta)h+g_n(\Delta))(2h-3)x-2f_n(\Delta)e4fx-\\
 & (f_n(\Delta)h+g_n(\Delta))4ey-2f_n(\Delta)e(-2h-3)y,
\end{align*}

so grouping all elements containing $y$ we get the following coefficient in fron of $y$
\begin{eqnarray*} 
4\Delta^ne+2f_n(\Delta)\Delta e-4g_n(\Delta)e+6f_n(\Delta)e+4f_n(\Delta)(-2e)=\\
2(2\Delta^n+f_n(\Delta)(\Delta-1)-2g_n(\Delta))e
 \end{eqnarray*} 
so $f_{n+1}(\Delta)=2\Delta^n+f_n(\Delta)(\Delta-1)-2g_n(\Delta).$

Grouping elements in front of $x$ we get
\begin{eqnarray*}
\Delta^n(2h-3)+f_n(\Delta)\Delta h+g_n(\Delta)\Delta-g_n(\Delta)2h+\\
3g_n(\Delta)-f_n(\Delta)h+4f_n(\Delta)h-2f_n(\Delta)h^2-2f_n(\Delta)4ef,
 \end{eqnarray*}
sum of last 3 terms is $-2f_n(\Delta)\Delta$, so
\begin{eqnarray*} 
g_{n+1}(\Delta)=-3\Delta^n+g_n(\Delta)\Delta+3g_n(\Delta)-2f_n(\Delta)\Delta=\\
-3\Delta^n+(\Delta +3)g_n(\Delta)-2f_n(\Delta)\Delta.
 \end{eqnarray*}
We easily see that $ f_n$ is a polynomial of degree $ n-1$ with the top coefficient $2n$, in particular any equation $ c=[z_1,x]+z_2 x$ has a unique solution in $z_1, z_2$ from $\mf{Z}(\mf{U}\mf{g}), $ provided that $c$ is of the form $2\psi ey+(h\psi +\psi_1)x$, where $\psi ,\psi_1$ are also from $\mf{Z}(\mf{U}\mf{g})$ (because $f_n$ form a basis of center, and by uniqueness in $z_1$ we mean uniqueness up to a constant), as an example of such $c$ we can take arbitrary linear combinations of elements of the form $\Delta^i[\Delta^j,x]$. Remark that only elements from $\mf{Z}(\mf{U}\mf{g})$ which commute with $x$ are constants. These observations will be needed later.
 
  Recall that we want to write $q_z$ as a commutator of $x$ with a central element of $\mf{U}\mf{g},$ let us rewrite $q_z$ as follows:
\begin{eqnarray*} 
z[\half\Delta,x]-(e[z,y]+\half h[z,x])+\half[z,x]=\\
z[\half\Delta,x]-[z,ey+\half hx]+\half[z,x]=z[\half\Delta,x]-[z,[\frac{1}{4}\Delta,x]]-[z,\frac{3}{4}x]+\\
\half[z,x]=\frac{1}{4}(z[\Delta,x]+[\Delta,x]z-[z,x]).
 \end{eqnarray*}
Since $[\Delta z,x]=[\Delta ,x]z+\Delta[z,x]$, its enough to show that there exists an element $z_1$ in the center of $\mf{U}\mf{g}$
so that $[z_1,x]=z[\Delta,x]-\Delta[z,x]=\Delta xz-zx\Delta$, let us take $z_1,z_2$ from center of $\mf{U}\mf{g}$ so that 
$z[\Delta,x]-\Delta[z,x]=[z_1,x]+z_2x$ (this is possible as remarked earlier). We want to show that $z_2$=0, applying $ad(f)$ to the last equality we get  $z[\Delta,y]-\Delta[z,y]=[z_1,y]+ z_2y$. Remark that $j(z[\Delta,x]-\Delta[z,x])=-(z[\Delta,y]-\Delta[z,y])$, therefore
\begin{eqnarray*}
[z_1,y]+z_2y=-j([z_1,x]+z_2x)=[z_1,y]-yz_2=[z_1+z_2,y]-z_2y.
 \end{eqnarray*}
Applying $ad(e)$ to the last line of equalities and using uniqueness of $z_2$ implies that $z_2=0$, and we have constructed a central element $t_z=t-\half hz-\omega_z$ of filtration degree 2. It remains to show that it generates the whole center. First claim is that
any element which commutes with $\mf{g}$ lies in the algebra generated by $t_z$ and $\mf{Z}(\mf{U}\mf{g}).$
 
Let   $a$ be such an element,  write $a$ as a polynomial in PBW basis (in $x,y$ writing $x$ on the right side of $y$ with coefficients in $\mf{U}\mf{g}$), and let $n$ be the smallest power of $x$ appearing in this polynomial, so we can write $a=bx^n$ since $ a$ commutes with $e$ and $x$ commutes with  $e$, implies that $b$ has weight $ -n$ and is maximal. However, since $H_z$ is union of finite dimensional modules over $ \mf{g}$ (recall that action is always adjoint here)  no maximal vectors can have negative weights, so  $n=0$. Let  $cy^m$ be a monomial of $a$ with highest power of $y$ appearing in monomials with no $x$ in it, then  $[e,a]=0$ implies that $[e,c]$=0 (becasue after taking commutator with $e$ powers of $x$ never decrease and power of $y$ never increase ), so $c$ is a maximal vector of $\mf{U}L$ of weight $m$. Now claim is that any maximal vector in $\mf{U}\mf{g}$ is generated by $e,\Delta$. Indeed, let $\alpha$ be such a vector, using induction (asuming homogenity) without loss of generality assume that it is not divisable by $e$, so if we write it in PBW basis (for $\mf{U}\mf{g}$) it will contain a monomial containing no $e$, so it has non positive wight so its weight is 0, therefore it lies in the center and we are done. So $m$ is even $m=2n$ and $c=\gamma e^{2n}$, $\gamma$ being central, now let us consider an  element $a-\gamma t_z^n$, call it $\zeta$. Then $\zeta$ commutes with $\mf{g}$ and has highest power of $y$ without $x$ less than $m$, now arguing inductively we will arrive at the central element of $\mf{U}\mf{g}$ so we are done. Finally, let $a$ be a central element of $H_z$, in particular, it commutes with $\mf{g}$ so it can be written as a polynomial of $ t_z$ with coefficients in $ \mf{Z}(\mf{U}\mf{g})$, let $\kappa t_z^n$ be the monomial of the top degree in $a$, since $[a,x]=0$, passing to the associated graded ring we get that $[x,\kappa]=0$, so $\kappa$ is a constant, so we  may disregard the top term, proceeding buy induction we are done.

\end{proof}

 To summerise, $t_z=t-\half hz-\omega_z$, where $t_z=ey^2+hxy-fx^2$, and $\omega_z$ is a polynomial in $\Delta$ with the property $[x,\omega_z]=\frac{1}{4}([\Delta z-z,x]+z[\Delta,x]-\Delta[z,x])$. In  what follows it will be convenient to use the following notation, denote by $F,G$ linear endomorphisms of $C[\Delta]$ such that $ F(\Delta^n)=f_n, G(\Delta^n)=g_n$. therefore $F,G$ are surjective linear endomorphisms  with kernel being constants. Thus in equation for $\omega_z$ after equating terms with $ey$ we get 
\begin{eqnarray*}
\frac{1}{4}(F(\Delta z-z)+zF(\Delta)-\Delta F(z))=-F(\omega_z), 
\end{eqnarray*}
now lets rewrite recursive relations on $f_n$ obtained before in terms of $F,G$, we get $F(\Delta z)=2z+(\Delta-1)F(z)-2G(z)$, so $\omega_z=-\frac{1}{4}F^{-1}(F(\Delta z-z)+zF(\Delta)-\Delta F(z))$, thus  
\begin{eqnarray*}
\omega_z &=&-\frac{1}{4}(2z+F(z)(\Delta-1)-2G(z)-F(z)+2z-\Delta F(z))\\
&=& -F^{-1}(z)+\half z+\half F^{-1}(G(z)), 
\end{eqnarray*}
as a consequence of this formula, top coefficient of $\omega_z$ is $-\frac{1}{2n}$ times the top coefficient of $z$ (assuming that the degree of $z$ is $n-1$ as a polynomial of $\Delta$), because according to computations of $f_n, g_n$, $F$ and $G$ decrease the degree of an input by 1, and top degree of $f_n$ is $2n$. 

   If $z$ is linear $z=a\Delta+b$, then explicit formula for a central element is the following:

\begin{eqnarray*}
 ey^2+hxy-fx^2-\half h(a\Delta+b)+\frac{1}{4}(a\Delta^2+(2b-a)\Delta).
 \end{eqnarray*} 
 \begin{theorem} 

   If $z$ is non-zero then $H_z$ has no outer derivations, and $H^1(H, H)$ is a rank 1 free module over the center of $H$.
\end{theorem}
\begin{proof}

  For the proof will need a description of maximal vectors of weight 1 in $H$. The claim is that this vector space as an algebra over the center of $H_z$ is generated by$[x,\mathbb{C}[\Delta]]$ and by $ \mathbb{C}[\Delta]x$. It suffices to show this when $z=0$, let $r=\sum \alpha_{ij}x^iy^j$ be a vector of maximal weight, we will use induction on total degree with respect to $x,y$ and biggest $j$ so that $\alpha_{0j}$ is not zero. If such coefficients are all $0$ then we can pull out $x$ we will be left with a maximal vector of weight 0, so in that case nothing is to prove. We see that $\alpha_{0j}$ is maximal, so  this term is equal to $\gamma e^ny^{2n-1}$, where $\gamma$ is from $\mathbb{C}[\Delta]$, so 
$r-\gamma t^{n-1}(\frac{1}{4}[\Delta,x])$ is also a maximal vector whose biggest power of $y$ without $x$ is lower than $j$, so continuing this procedure all is left to show is that $\alpha[\Delta,x]$ can be expressed appropriately, where $\alpha$ is a polynomial of $\Delta$. Indeed, this element can be clearly written as $[\Delta,\gamma]+\omega x$, where$ \omega, \gamma$ are polynomials of $\Delta$ (as explained earlier), so we are done.

 We claim that $H^1(H, H)$ is rank 1 free module over the center, indeed assume that $D$ is a derivation of $H$, we may assume that it vanishes on $\mf{U}\mf{g}$ (consequence of simplicity of the Lie algebra $\mf{g}$), so $D$ is a $\mf{g}$-map, therefore $Dx$ is a maximal vector of rank 1, so by modifying with inner derivation, we may assume that $Dx=\sum t^i \alpha_i x$, where  is  $\alpha_i$ is a polynomial of $\Delta$. We claim that $\alpha_i$ is constant, indeed, assume the contrary, then

\begin{eqnarray*}
 D(xy)=\sum (t^i\alpha_ixy+xt^i\alpha_iy)=\sum t^i(2\alpha_ixy-[\alpha_i,x]y)\\
 D(yx)=\sum (t^i\alpha_iyx+yt^i\alpha_ix)=\sum t^i(2\alpha_ixy-[\alpha_i,y]x).
\end{eqnarray*} 
 Thus $\sum t^i[\alpha_i,x]y=\sum t^i[\alpha_i,y]x$, but this can not happen becaus after equating highest powers of $t$, left hand side contains terms with $y^2$ and right hand side does not, this implies that all $\alpha_i$ are constants, so $D$ is given by $\mf{g}$-invariant Euler derivative  on  $\mf{U}(\mf{g}\ltimes V)$ (with respect to $x,y$) multiplied by a central element, therefore we are done.
  Now lets consider remaining cases ($z$ non zero). We claim that every derivation is inner, suppose $D$ is a derivation of $H$, again we may assume that it vanishes on $\mf{U}\mf{g}$. As before, we may assume that $ Dx=\sum t_z^i\alpha_ix $, where $\alpha_i$  are central in $\mf{U}\mf{g}$ and  $\alpha_0$ is non zero, so $Dy=\sum t_z^i\alpha_iy$. Thus applying $D$ to $z$ we get

$0=Dz=[Dx,y]+[x,Dy]=2\sum t_z^i\alpha_iz+\sum t_z^i([\alpha_i,y]x+[x,\alpha_i]y)$
 so 
\begin{eqnarray*}
 2\sum t_z^i\alpha_iz=\sum t_z^i([\alpha_i,x]y-[\alpha_i,y]x)
\end{eqnarray*}
For any $\alpha$ from $\C[\Delta]$, we have
\begin{eqnarray*} 
   && [\alpha,x]y-[\alpha,y]x \\
    &=& (2F(\alpha)ey+(F(\alpha)h+G(\alpha))x)y-(2F(\alpha)fx+(-F(\alpha)h+G(\alpha))y)x \\
    &=& 2F(\alpha)(ey^2+hxy-fx^2-\half hz)+G(\alpha)z=2F(\alpha)(t_z+\omega_z)+G(\alpha)z 
\end{eqnarray*}
 so our equality is as follows: $ 2\sum t_z^i\alpha_iz=\sum t_z^i(2F(\alpha_i)(t+\omega_z)+G(\alpha_i)z)$, now lets equate terms with no $t_z$, we get $2\alpha_0z=2F(\alpha_0)B(z)+G(\alpha_0)z$. We have an equality in polynomials in $\Delta$, so $\alpha_0$ is not a constant otherwise right hand side would be 0, now comparing coefficients of top degrees in $\Delta$ gives a contradiction, because as we have observed $ F$ multiplies the top coefficient of an input by a positive number, however the top coefficient of $\omega_z$ is equal to a negative number times the top coefficient of $z$. End of the proof.

\end{proof}

        The University of Chicago, Department of Mathematics, Chicago, Il 60637, USA.

\end{document}